\documentclass[a4paper,12pt]{amsart}
\usepackage{a4,amsmath,amssymb,amscd}

\theoremstyle{plain}

\newtheorem{theo}{Theorem}[section]

\newtheorem{lemm}[theo]{Lemma}
\newtheorem{coro}[theo]{Corollary}

\theoremstyle{definition}
\newtheorem*{defi}{Definition}
\newtheorem{exam}[theo]{Example}

\theoremstyle{remark}
\newtheorem*{rema}{Remark}

\numberwithin{equation}{section}

\newcommand{\field}[1]{\mathbb{#1}}
\newcommand{\C}{\field{C}}
\newcommand{\R}{\field{R}}
\newcommand{\Z}{\field{Z}}
\newcommand{\Q}{\field{Q}}
\newcommand{\I}{\mathcal I}
\renewcommand{\k}{\mathbf k}
\newcommand{\w}{\omega}
\newcommand{\oP}{\overline P}
\newcommand{\A}{\mathcal A}
\newcommand{\oA}{\bar{A}}
\newcommand{\ly}{\lambda(y,y')}

\DeclareMathOperator{\rank}{rank}
\DeclareMathOperator{\ind}{Ind}
\newcommand{\indT}{\ind_T}
\DeclareMathOperator{\sgn}{\epsilon}

\begin{document}

\title{$h$-vectors of Gorenstein* simplicial posets}

\author[M.~Masuda]{Mikiya Masuda}
\address{Department of Mathematics, Graduate School of Science, 
Osaka City University, Sugimoto, Sumiyoshi-ku, Osaka 558-8585, Japan}
\email{masuda@sci.osaka-cu.ac.jp}

\maketitle

\section{Introduction}
A \emph{simplicial poset} $P$ 
(also called a \emph{boolean poset} and a \emph{poset of boolean type}) 
is a finite poset with a smallest element $\hat 0$ such that every 
interval $[\hat 0,y]$ for $y\in P$ is a boolean algebra, i.e., $[\hat 0,y]$ 
is isomorphic to the set of all subsets of a finite set, ordered by 
inclusion.  
The set of all faces of a (finite) simplicial complex 
with empty set added forms a simplicial poset ordered by inclusion, 
where the empty set is the smallest element.  
Such a simplicial poset is called the \emph{face poset} of a 
simplicial complex, and two simplicial complexes are isomorphic if and 
only if their face posets are isomorphic.  Therefore, a simplicial 
poset can be thought of as a generalization of a simplicial complex. 

Although a simplicial poset $P$ is not necessarily the face poset of 
a simplicial complex, it is always the face poset of a CW-complex 
$\Gamma(P)$. 
In fact, to each $y\in P\backslash\{\hat 0\}=\oP$, 
we assign a (geometrical) simplex whose face poset is $[\hat 0,y]$ and 
glue those geometrical simplices according to the order relation in $P$.
Then we get the CW-complex $\Gamma(P)$ such that 
all the attaching maps are inclusions.  
For instance, if two simplices of a same dimension are identified on 
their boundaries via the identity map, then it is not a simplicial complex 
but a CW-complex obtained from a simplicial poset. 
The CW-complex $\Gamma(P)$ has a well-defined barycentric subdivision 
which is isomorphic to 
the order complex $\Delta(\oP)$ of the poset $\oP$.  Here 
$\Delta(\oP)$ is a simplicial 
complex on the vertex set $\oP$ whose faces are the chains of $\oP$. 

We say that 
$y\in P$ has rank $i$ if the interval $[\hat 0,y]$ 
is isomorphic to the boolean algebra of rank $i$ (in other words, the 
face poset of an $(i-1)$-simplex), and the rank of $P$ is 
defined to be the maximum of ranks of all elements in $P$. 
Let $d=\rank P$. In exact analogy to simplicial complexes,  
the $f$-vector of the simplicial poset 
$P$, $(f_0,f_1,\dots,f_{d-1})$, is defined by 
\[
f_i=f_i(P)=\#\{ y\in P \mid \rank y=i+1\}
\]
and the $h$-vector of $P$, $(h_0,h_1,\dots,h_d)$, 
is defined by the following identity:
\begin{equation*} \label{eqn:fh}
\sum_{i=0}^d f_{i-1}(t-1)^{d-i}=\sum_{i=0}^d h_i t^{d-i},
\end{equation*}
where $f_{-1}=1$, so $h_0=1$. 
Note that the number of facets of $P$, that is $f_{d-1}$, is related to 
$h$-vectors as follows: 
\begin{equation} \label{eqn:fn-1}
f_{d-1}=\sum_{i=0}^d{h_i}.
\end{equation}
When $P$ is the face poset of a simplicial complex $\Sigma$, the 
$f$- and $h$-vector of $P$ coincide with the classical $f$- and $h$-vector of 
the simplicial complex $\Sigma$ respectively.  

$f$-vectors and $h$-vectors have 
equivalent information, but $h$-vectors are often easier than $f$-vectors. 
In \cite{stan91}, R.~Stanley discussed characterization of $h$-vectors 
for certain classes of simplicial posets.  
For example, he proved that 
a vector $(h_0,h_1,\dots,h_d)$ of integers with $h_0=1$ is the $h$-vector of 
a Cohen-Macaulay simplicial poset of rank $d$ 
if and only if $h_i\ge 0$ for all $i$.  
Gorenstein* simplicial posets are more special than Cohen-Macaulay 
simplicial posets.  
If the CW-complex 
$\Gamma(P)$ is homeomorphic to a sphere of dimension $d-1$, then 
the simplicial poset $P$ of rank $d$ is Gorenstein* 
(see Section~\ref{sect:gore} for more details). 
It is known that $h$-vectors of Gorenstein* simplicial posets 
satisfy Dehn-Sommerville equations 
$h_i=h_{d-i}$ for all $i$,
in addition to the non-negativity conditions $h_i\ge 0$.  
In this paper we will prove that $h$-vectors of Gorenstein* simplicial 
posets must satisfy one more subtle 
condition conjectured by Stanley in \cite{stan91}, see 
\cite{duva94}, \cite{ma-pa03}, \cite{stan91} for partial results. 

\begin{theo} \label{theo:main}
If $P$ is a Gorenstein* simplicial poset of rank $d$ and $h_i(P)=0$ 
for some $i$ between $0$ and $d$, then $\sum_{i=0}^d h_i(P)$, that 
is the number of facets of $P$ by (\ref{eqn:fn-1}), is even. 
\end{theo}

Combining this with Theorem 4.3 in \cite{stan91}, one completes 
characterization of $h$-vectors of Gorenstein* simplicial posets. 

\begin{coro} \label{coro:characterization}
Let $(h_0,h_1,\dots,h_d)$ be a vector of non-negative integers with 
$h_i=h_{d-i}$ for all $i$ and $h_0=1$.  There is 
a Gorenstein* simplicial poset $P$ of rank $d$ with 
$h_i(P)=h_i$ for all $i$ if and only if either $h_i>0$ for all $i$, or else 
$\sum_{i=0}^d h_i$ is even. 
\end{coro}

Our proof of Theorem~\ref{theo:main} is purely algebraic but the idea 
stems from topology, so we will explain how our proof is related 
to topology in Section~\ref{sect:topo}.  
A main tool to study the $h$-vector of a simplicial poset $P$ is 
a (generalized) face ring $A_P$ introduced in \cite{stan91} of the poset 
$P$.  
In Section~\ref{sect:rest} we discuss restriction maps 
from $A_P$ to polynomial rings. In Section~\ref{sect:index} we construct 
a map called an index map from $A_P$ to a polynomial ring.  
Theorem~\ref{theo:main} is proven in Section~\ref{sect:gore}.  

I am grateful to Takayuki Hibi for informing me of the above problem 
and for his interest.  I am also grateful to Akio Hattori and 
Taras Panov for the successful collaborations (\cite{HM}, \cite{ma-pa03}) 
from which the idea used in this paper originates.  
Finally I am grateful to Ezra Miller for his comments on an eariler 
version of the paper, which were very helpful to improve the paper.

\section{Relation to topology} \label{sect:topo}

In the toric geometry, simplicial convex polytopes 
are closely related to \emph{toric} manifolds or 
orbifolds (see \cite{fult93}).  Similarly to this, 
Gorenstein* simplicial posets, which contain the boundary complexes of 
simplicial polytopes as examples, are closely related to 
objects (in topology) called \emph{torus} manifolds 
or orbifolds (see \cite{HM}, \cite{ma-pa03}), and the proof 
of Theorem~\ref{theo:main} is motivated by a topological observation 
described in this section.  
Here a torus manifold (resp. orbifold) means a 
closed smooth manifold (resp. orbifold) of dimension $2d$ with an effective 
smooth action of a $d$-dimensional torus group having at least one fixed 
point. 

We shall illustrate relations between combinatorics 
and topology with simple examples. 
In the following, $T$ will denote 
the product of $d$ copies of the circle group consisting of complex 
numbers with unit length, i.e., $T$ is a $d$-dimensional torus group. 

\begin{exam} \label{CPd}
A complex projective space $\C P^d$ has a $T$-action 
defined in the homogeneous coordinates by
\[
(t_1,\dots,t_d)\cdot(z_0:z_1:\dots:z_d)=(z_0:t_1z_1:\dots:t_dz_d).
\]
The orbit space $\C P^d/T$ has a natural face structure.  Its facets 
are the images of (real) codimension two submanifolds $z_i=0$ 
$(i=0,1,\dots,d)$ under the quotient map $\C P^d\to \C P^d/T$. 
The map (called a moment map)
\[
(z_0:z_1:\dots:z_d) \mapsto 
\frac{1}{\sum_{i=0}^d |z_i|^2}(|z_1|^2,\dots,|z_d|^2)
\]
induces a face preserving homeomorphism from the orbit space $\C P^d/T$ 
to a standard $d$-simplex. 
The face poset of $\C P^d/T$ ordered by reverse inclusion 
(so $\C P^d/T$ iteself is the smallest element) 
is the face poset of a simplicial complex of dimension $d-1$ and 
Gorenstein*.  

Similarly, the product of $d$ copies of $\C P^1$ 
admits a $T$-action, the orbit space $(\C P^1)^d/T$ is 
homeomorphic to a $d$-cube, and 
the face poset of $(\C P^1)^d/T$ ordered by reverse inclusion 
is also the face poset of  
a simplicial complex of dimension $d-1$ and Gorenstein*.  

In any case, the orbit space is a simple convex polytope and its polar 
is a simplicial convex polytope.  
The Gorenstein* simplicial complex is the \emph{boundary} 
complex of the simplicial convex polytope. 
\end{exam}

\begin{exam}\label{2dsphere}
Let $S^{2d}$ be the $2d$-sphere  
identified with the following subset in $\C^{d}\times\R$:
$$
  \bigl\{ (z_1,\dots,z_d,y)\in\C^d\times\R\mid |z_1|^2+\dots+|z_d|^2+y^2=1 
  \bigr\},
$$
and define a $T$-action on $S^{2d}$ by
$$
  (t_1,\dots,t_d)\cdot(z_1,\dots,z_d,y)=(t_1z_1,\dots,t_dz_d,y).
$$
The facets in the orbit space $S^{2d}/T$ are the images of codimension 
two submanifolds $z_i=0$ $(i=1,\dots,d)$ 
under the quotient map $S^{2d}\to S^{2d}/T$, and the map
\[
(z_1,\dots,z_d,y)\to (|z_1|,\dots,|z_d|,y)
\]
induces a face preserving homeomorphism from 
$S^{2d}/T$ to the following subset of the $d$-sphere: 
\[
\{ (x_1,\dots,x_d,y)\in \R^{d+1}\mid x_1^2+\dots+x_d^2+y^2=1,\ x_1\ge 0,
\dots,x_d\ge 0\}.
\]
The orbit space $S^{2d}/T$ is not (isomorphic to) 
a simple convex polytope  because 
the intersection of $d$ facets consists of two points, 
but it is a manifold with corners 
and every face (even $S^{2d}/T$ itself) is acyclic.  
The face poset of $S^{2d}/T$ ordered by reverse inclusion 
is not the face poset of a simplicial complex when $d\ge 2$. 
However, it is a simplicial poset and Gorenstein*.  The geometric 
realization of the face poset of $S^{2d}/T$ 
is formed from two $(d-1)$-simplices by gluing their boundaries via 
the identity map.  

\end{exam}

More generally, it is proved in \cite{ma-pa03} that 
if a torus manifold $M$ has vanishing odd degree cohomology, 
then the orbit space $M/T$ is a manifold with corners 
and every face (even $M/T$ itself) is acyclic; so the face poset 
of $M/T$ ordered by reverse inclusion is a Gorenstein* simplicial 
poset, say $P$.  Moreover, $h_i(P)$ agrees with the $2i$-th betti number 
$b_{2i}(M)$ of $M$ and the equivariant cohomology ring 
$H^*_T(M;\Z)$ of $M$ is isomorphic to the face ring $A_P$ of $P$ 
(defined over $\Z$). 
Here $H^*_T(M;\k)$ for a ring $\k$ is defined as 
\[
H^*_T(M;\k):=H^*(ET\times_T M;\k)
\]
where $ET$ is the total space of the universal principal $T$-bundle 
(on which $T$ acts freely) and $ET\times_T M$ is the orbit space of 
the product $ET\times M$ by the diagonal $T$-action.

A projective toric orbifold is related to a simplicial convex polytope 
as in Example~\ref{CPd}, and the $h$-vector of the simplicial convex polytope 
agrees with the (even degree) betti numbers of 
the toric orbifold.  Noting this fact, 
Stanley \cite{stan80} deduced constraints on the $h$-vector 
by applying the hard Lefschetz theorem to 
the toric orbifold and completed the characterization of 
$h$-vectors of simplicial convex polytopes.  
In some sense our proof of Theorem~\ref{theo:main} is on this line.  
The topological argument developed below in this section 
is not complete but would be helpful 
for the reader to understand what is done in subsequent sections. 

Let $P$ be a Gorenstein* simplicial poset of dimension $d-1$.  
Looking at the result in \cite{ma-pa03} mentioned above, it is likely 
that there exists a torus orbifold $M$ which have the following properties: 

\medskip
\noindent
{\bf Properties.} 
\begin{enumerate}
\item $H^{odd}(M;\Q)=0$, 
\item $h_i(P)=b_{2i}(M)$, 
\item $H^*_T(M;\Q)$ is isomorphic to 
$A_P$ (defined over $\Q$). 
\end{enumerate}

\medskip

What we will use to deduce the necessity in Theorem~\ref{theo:main} 
is the \emph{index map} (or evaluation map) in equivariant cohomology: 
\[
\indT\colon H^*_T(M;\Q)\to H^{*-2d}_T(pt;\Q)=H^{*-2d}(BT;\Q),
\]
where $BT=ET/T$ is the classifying space of principal $T$-bundles.  
The index map is nothing but the Gysin homomorphism in 
equivariant cohomology induced from the collapsing map $\pi\colon M\to pt$. 
As is well-known, $BT$ is the product of $d$ copies of $\C P^\infty$ 
(up to homotopy) 
and $H^*(BT;\Q)$ is a polynomial ring in $d$ variables of degree two. 
The index map $\indT$ 
decreases cohomological degrees by $2d$ because the dimension of $M$ is 
$2d$. Moreover, $H^*_T(M;\Q)$ is a module over $H^*(BT;\Q)$ through 
$\pi^*\colon H^*(BT;\Q)=H^*_T(pt;\Q)\to H^*_T(M;\Q)$ and 
$\indT$ is an $H^*(BT;\Q)$-module map.  Since $H^{odd}(M;\Q)=0$ and 
$H^*(BT;\Q)$ 
is a polynomial ring in $d$ variables, say $t_1,\dots,t_d$, 
the quotient ring of $H^*_T(M;\Q)$ by the ideal generated 
by $\pi^*(t_1), \dots, \pi^*(t_d)$ agrees with the ordinary cohomology 
$H^*(M;\Q)$.  Similarly, the quotient ring of $H^*_T(pt;\Q)=H^*(BT;\Q)$ 
by the ideal 
generated by $t_1,\dots,t_d$ agrees with $H^*(pt;\Q)$.  Therefore 
the index map in equivariant cohomology 
induces the index map in ordinary cohomology:
\[
\ind \colon H^*(M;\Q) \to H^{*-2d}(pt;\Q). 
\]
This map agrees with the Gysin homomorphism in ordinary cohomology 
induced from the collapsing map $\pi$, so it is 
the evaluation map on a fundamental class of $M$. 
Thus, we have a commutative diagram:
\[
\begin{CD}
H^{2d}_T(M;\Q)@>\indT>>H^{0}_T(pt;\Q)=H^{0}(BT;\Q)=\Q\\
@VVV          @VVV \\
H^{2d}(M;\Q)@>\ind >> H^{0}(pt;\Q)=\Q,
\end{CD}
\]
where the right vertical map is the identity. 

A key thing is to find an element $\w_T$ in $H^{2d}_T(M;\Q)$ such that 
\begin{enumerate}
\item[(i)] $\w_T$ is a polynomial in elements of $H^2_T(M;\Q)$,
\item[(ii)] 
$\indT(\w_T)$ is an integer and $\indT(\w_T)\equiv \chi(M) \pmod 2$, 
where $\chi(M)$ is the Euler characteristic of $M$.  
\end{enumerate}
We may think of $\w_T$ as a \lq\lq lifting" of the equivariant 
top Stiefel-Whitney class $w_{2d}^T(M)\in H^{2d}_T(M;\Z/2)$ of $M$.  
If we find such an element $\w_T$, then it follows from 
the commutativity of the above diagram that  
\begin{equation} \label{eqn:CD}
\indT(\w_T)=\ind(\w)
\end{equation}
where $\w$ is the image of $\w_T$ 
under the left vertical map in the above diagram. 

Now suppose $h_i(P)=0$ for some $1\le i\le d-1$. 
Then the $2i$-th betti number $b_{2i}(M)$ of $M$ is zero 
by property (2) and the element $\w$ vanishes 
because it is a polynomial in degree two elements by (i) above, 
so the right hand side 
of (\ref{eqn:CD}) is zero and $\chi(M)$ is even by (ii) above. 
On the other hand, it follows from properties (1) and (2) that 
$$\chi(M)=\sum_{i=0}^d b_{2i}(M)=\sum_{i=0}^d h_i(P).$$ 
These prove that $\sum_{i=0}^d h_i(P)$ is even. 

It turns out that the argument developed above works without assuming 
the existence of the torus orbifold $M$.  
In fact, the face ring $A_P$ takes the place 
of $H^*_T(M;\Q)$ by property (3) 
and an l.s.o.p. for $A_P$ plays the role of  
$\pi^*(t_1),\dots,\pi^*(t_d)$ so that the polynomial ring generated 
by the l.s.o.p. corresponds to the polynomial ring $\pi^*(H^*(BT;\Q))$ 
(or $H^*(BT;\Q)$ since $\pi^*$ is injective). 
The index map $\indT$ has an expression (so-called Lefschetz fixed point 
formula) 
in terms of local data around $T$-fixed 
points of $M$, and since the formula is purely algebraic, 
one can use it to define an \lq\lq index map" from $A_P$.  
To carry out this idea, we need to study restriction maps from $A_P$ 
to polynomial rings because restriction maps to $T$-fixed points 
in equivariant cohomology are involved in the Lefschetz fixed point formula.  
We will discuss such restriction maps in Section~\ref{sect:rest} and 
construct the index map from $A_P$ in Section~\ref{sect:index}.

\section{Restriction maps} \label{sect:rest}

In this and next sections, we consider rings over $\Q$. 
A main tool to study the $h$-vector of a (finite) simplicial poset 
$P$ is the face ring 
$A_P$ of the poset $P$ introduced by Stanley in \cite{stan91}. 
We recall it first. 

\begin{defi} 
Let $P$ be a simplcial poset of rank $d$ 
with elements $\hat 0, y_1,\dots,y_p$. 
Let $\oA=\Q[y_1,\dots,y_p]$ be the polynomial ring over $\Q$ in the 
variables $y_i$ and define $\I_P$ to be the ideal of $\oA$ generated by 
the following elements:
$$y_iy_j-(y_i\wedge y_j)(\sum_z z),$$
where $y_i\wedge y_j$ is 
the greatest lower bound of $y_i$ and $y_j$, $z$ ranges over all 
minimal upper bounds of $y_i$ and $y_j$, and we understand $\sum_z z=0$ 
if $y_i$ and $y_j$ have no common upper bound.  Then the 
face ring $A_P$ of the simplicial poset $P$ 
is defined as the quotient ring $\oA/\I_P$ and made graded 
\[
A_P=(A_P)_0\oplus (A_P)_1\oplus \dots
\]
by defining $\deg y_i=\rank y_i$.  
The ring $A_P$ reduces to a classical Stanley-Reisner face ring 
when $P$ is the face poset of a simplicial complex. 
\end{defi}

We denote by $P_s$ the subset of $P$ consisting of elements of $\rank s$.
Elements in $P_1$ will be denoted by $x_1,\dots,x_n$ and called 
\emph{atoms} in $P$. The set $\{x_1,\dots,x_n\}$ is a basis of $(A_P)_1$. 

Suppose that $y$ is an element of $P_d$. Then the interval $[\hat 0,y]$ is a 
boolean algebra of rank $d$ and 
$A_{[\hat 0,y]}$ is a polynomial ring in $d$ variablesD
Sending all elements in $P$ which are not lower than $y$ to zero, 
we obtain an epimorphism 
\[
\iota_y\colon A_P\to A_{[\hat 0,y]}.
\]
Since $\Q$ is a field with infinitely many elements, $A_P$ admits an l.s.o.p. 
$\theta_1,\dots,\theta_d$ (see the proof of Theorem 3.10 in \cite{stan91}). 
In the following we fix the l.s.o.p. and denote by $\Theta$ the vector space 
of dimension $d$ spanned by $\theta_1,\dots,\theta_d$ over $\Q$, and by 
$\Q[\Theta]$ the polynomial ring generated by $\theta_1,\dots,\theta_d$. 
Note that $\Theta$ is a vector subspace of $(A_P)_1$ and $\Q[\Theta]$ 
is a subring of $A_P$. 

\begin{lemm} \label{lemm:iso}
The restriction of $\iota_y$ to 
$\Q[\Theta]$ is an isomorphism onto $A_{[\hat 0,y]}$. 
\end{lemm}

\begin{proof} 
Since $A_P$ is finitely generated as a 
$\Q[\Theta]$-module, so is $A_{[\hat 0,y]}$. 
This implies that $\iota_y$ maps the vector space 
$\Theta$ isomorphically onto the vector 
space spanned by $d$ elements of degree one generating the polynomial ring 
$A_{[\hat 0,y]}$, thus the lemma follows. 
\end{proof}

Henceforth, we identify 
$A_{[\hat 0,y]}$ with $\Q[\Theta]$ via 
$\iota_y$, and think of $\iota_y$ as a map to 
$\Q[\Theta]$, i.e., 
\[
\iota_y\colon A_P\to \Q[\Theta].
\]  
Note that $\iota_y$ is the identity on the subring $\Q[\Theta]$ and 
a $\Q[\Theta]$-module mapD

For $w\in P_s$, we set 
\[
\A(w):=\{ i\in \{1,\dots,n\}\mid \text{$x_i$ is an atom lower than $w$}\}.
\]
The cardinality of $\A(w)$ is $s$. 
Let $y\in P_d$.  By definition of $\iota_y$, 
\begin{equation} \label{eqn:iota}
\text{$\iota_y(x_i)=0$ whenever 
$i\notin \A(y)$.}
\end{equation}  
We set 
\begin{equation} \label{eqn:theta}
\theta_i(y):=\iota_y(x_i)\quad\text{for $i\in \A(y).$}
\end{equation}
Since $\iota_y\colon (A_P)_1\to \Theta$ is surjective and the cardinality 
of $\A(y)$, that is $d$, agrees with the dimension of $\Theta$, the set 
$\{\theta_i(y)\mid i\in \A(y)\}$ is a basis of $\Theta$. 

Let $z\in P_{d-1}$.  Let $y$ be an element in $P_d$ above $z$ and 
define $\ell\in\{ 1,\dots,n\}$ by 
\[
\A(y)\backslash \A(z)=\{ \ell\}.
\]
The canonical map $A_{[\hat 0,y]}=\Q[\Theta] \to A_{[\hat 0,z]}$ is surjective 
and $A_{[\hat 0,z]}$ can canonically 
be identified with $\Q[\Theta]/(\theta_\ell(y))$. 
Let $y'$ be another element in $P_d$ above $z$ and define 
$\ell'\in\{ 1,\dots,n\}$ similarly to $\ell$.  
It may happen that $\ell=\ell'$.  Since 
\begin{equation} \label{eqn:Az}
\Q[\Theta]/(\theta_\ell(y))=A_{[\hat 0,z]}=\Q[\Theta]/(\theta_{\ell'}(y')), 
\end{equation}
$\theta_\ell(y)$ and $\theta_{\ell'}(y')$ are same up to a non-zero scalar 
multiple; so the following lemma makes sense. 

\begin{lemm} \label{lemm:5}
$\iota_y(\alpha)\equiv \iota_{y'}(\alpha) \mod{\theta_{\ell}(y)}$ 
for any $\alpha\in A_P$. In particular, $\theta_i(y)\equiv \theta_i(y') 
\mod{\theta_\ell(y)}$ for $i\in \A(z)(=\A(y)\backslash\{\ell\}=\A(y')
\backslash\{\ell'\})$. 
\end{lemm} 

\begin{proof}
We have canonical surjections $A_P\to A_{[\hat 0,y]}\to A_{[\hat 0,z]}$ and 
$A_P\to A_{[\hat 0,y']}\to A_{[\hat 0,z]}$, whose composite surjections 
$A_P\to A_{[\hat 0,z]}$ are the same.  Therefore the lemma follows from 
(\ref{eqn:Az}). 
\end{proof}

\bigskip
\section{Index maps} \label{sect:index}

In this section, we define an \lq\lq index map" 
from $A_P$ to the polynomial 
ring $\Q[\Theta]$, which corresponds to the  
index map $\indT$ in Section~\ref{sect:topo}. 
It is a $\Q[\Theta]$-module map, so  
it induces a homomorphism from the quotient 
$A_P/(\Theta)$ modulo the linear system of parameters 
$\theta_1,\dots,\theta_d$ to $\Q$.  This induced map 
corresponds to the index map $\ind$ in Section~\ref{sect:topo}. 

We shall make some observations needed later before we define the index map. 
Let $z\in P_{d-1}$ and let $y,y'\in P_d$ lie above $z$ as before.  
Give an orientation on $\Theta$ 
determined by an ordered basis $(\theta_1,\dots,\theta_d)$ and 
choose an order of the basis $\{ \theta_i(y) \mid i\in \A(y)\}$ whose 
induced orientation on $\Theta$ agrees with the given orientation. 
We then define $m(y)$ to be 
the determinant of a matrix sending the ordered 
basis $\{\theta_i(y)\mid i\in \A(y)\}$
to the ordered basis $(\theta_1,\dots,\theta_d)$.  Note that $m(y)$ 
is positive. 
It follows from the latter statement in Lemma~\ref{lemm:5} that 
\begin{equation} \label{eqn:|b|'}
m(y)\theta_\ell(y)=\ly m(y')\theta_{\ell'}(y'), 
\end{equation}
where $\ly=\pm 1$.  
If $\A(y)=\A(y')$, then $\ell=\ell'$ and both $\theta_\ell(y)$ and 
$\theta_{\ell'}(y')$ restrict to the element $x_\ell$ 
in $A_{[\hat 0,x_\ell]}$.  Therefore 
\begin{equation} \label{eqn:m(y)=m(y')}
\text{$m(y)=m(y')$\ \ (and $\ly=1$) \ if\ \ $\A(y)=\A(y')$.}
\end{equation}

The order of the basis $\{ \theta_i(y) \mid i\in \A(y)\}$ 
determines an order of atoms $x_i$ $(i\in \A(y))$ and then determines 
an orientation on the $(d-1)$-simplex with those atoms as vertices.  
The oriented $(d-1)$-simplex obtained in this way 
is denoted by $\langle y\rangle$.  Then the boundaries  
$\partial \langle y\rangle$ and $\partial \langle y'\rangle$ 
of $\langle y\rangle$ and $\langle y'\rangle$ have 
opposite orientations on the $(d-2)$-simplex $[z]$ corresponding to $z$ 
(in other words, $[z]$ does not appear in 
$\partial\langle y\rangle+\partial \langle y'\rangle$) 
if and only if $\ly=-1$. 

Now we pose the following assumption, which we shall see in 
Section~\ref{sect:gore} is satisfied by all Gorenstein* simplicial posets. 

\medskip
\noindent
{\bf Assumption.}
\begin{enumerate}
\item For any $z\in P_{d-1}$, there are exactly two elements in 
$P_d$ above $z$. 
\item One can assign a sign $\sgn(y)\in \{ \pm 1\}$
to each $y\in P_d$ so that 
$
\sum_{y\in P_d}\sgn(y)\langle y\rangle 
$
is a cycle (hence defines a fundamental class in $H_{d-1}(\Gamma(P);\Z)$
where $\Gamma(P)$ denotes the CW-complex explained in the Introduction).  
\end{enumerate}

\medskip
\noindent
When $\langle y\rangle$ and $\langle y'\rangle$ share 
a $(d-2)$-simplex $[z]$, it follows from the above assumption that 
$[z]$ does not appear in 
$\partial(\sgn(y)\langle y\rangle)+\partial(\sgn(y')\langle y'\rangle)$. 
Therefore, 
\begin{equation} \label{eqn:sign of b}
\text{$\ly$ and $\sgn(y)\sgn(y')$ have opposite signs}
\end{equation}
by the remark mentioned above the Assumption.

\begin{defi} 
For a simplicial poset $P$ which satisfies the Assumption above, 
we define the \emph{index map} by 
\begin{equation} \label{eqn:index}
\indT(\alpha):= \sum_{y\in P_d}\frac{\sgn(y)\iota_y(\alpha)}
{m(y)\prod_{i\in \A(y)}\theta_i(y)}
\qquad\text{for $\alpha\in A_P$.}
\end{equation}
\end{defi}

Apparently, $\indT(\alpha)$ lies in the quotient field of 
$\Q[\Theta]$, but we have 

\begin{theo} 
$\indT(\alpha)\in \Q[\Theta]$ for any $\alpha\in A_P$. 
\end{theo}

\begin{rema}
The proof given below is essentially same as that of Theorem 2.2 in \cite{GZ}. 
A similar result can be found in \cite[Section 8]{HM}. 
\end{rema}

\begin{proof} 
The right hand side of (\ref{eqn:index}) can be expressed as 
\begin{equation} \label{eqn:common}
\frac{g}{\prod_{j=1}^N f_j}
\end{equation}
with $g\in \Q[\Theta]$ and 
$f_j\in \Theta\subset\Q[\Theta]$ 
such that any two of $f_1,\dots,f_N$ are linearly independent. 
It suffices to show that $f_1$ divides $g$. 

Let $Q$ be the set of $y\in P_d$ such that 
$\theta_i(y)$ is not a scalar multiple of $f_1$ for every $i\in \A(y)$, 
and let $Q^c$ be the complement of $Q$ in $P_d$. 
In (\ref{eqn:index}), the sum of terms for elements in $Q$ reduces to 
\begin{equation} \label{eqn:Q}
\sum_{y\in Q}\frac{\sgn(y)\iota_y(\alpha)}{m(y)\prod_{i\in \A(y)}
\theta_i(y)}=\frac{g_1}{\prod_{j=2}^N f_j} 
\end{equation}
with $g_1\in \Q[\Theta]$, 
so that $f_1$ does not appear in the denominator. 

On the other hand, if $y\in Q^c$, then it follows from the 
definition of $Q$ that there is an element $\ell\in \A(y)$ such that 
\begin{equation} \label{eqn:c}
\theta_\ell(y)=cf_1 \qquad(0\not=c\in \Q),
\end{equation}
and there is a unique element 
$z\in P_{d-1}$ such that $z$ is lower than $y$ and 
$\A(z)=\A(y)\backslash \{\ell\}$. By assumption, there is a unique 
element in $P_d$ which lies above $z$ and is different from $y$. 
We denote it by $y'$. 
Now we are in the same situation as before. 
It follows from (\ref{eqn:|b|'}) and (\ref{eqn:c}) that 
$y'$ is also an element in $Q^c$. 
Noting that $\A(y)=\A(z)\cup\{\ell\}$ and $\A(y')=\A(z)
\cup \{\ell' \}$ and using (\ref{eqn:|b|'}), 
we combine the two terms in (\ref{eqn:index}) for $y$ and $y'$ to get 
\begin{equation} \label{eqn:y,y'}
\begin{split}
&\frac{\sgn(y)\iota_y(\alpha)}{m(y)\prod_{i\in \A(y)}\theta_i(y)}+
\frac{\sgn(y')\iota_{y'}(\alpha)}{m(y')\prod_{i\in \A(y')}\theta_i(y')}\\
=&\frac{\sgn(y)\iota_y(\alpha)\prod_{i\in \A(z)}
\theta_i(y')+\ly\sgn(y')\iota_{y'}(\alpha)
\prod_{i\in \A(z)}\theta_i(y)}
{m(y)\theta_\ell(y)\prod_{i\in \A(z)}\theta_i(y)
\prod_{i\in \A(z)}\theta_i(y')}.
\end{split}
\end{equation}
Here 
\[
\iota_y(\alpha)\prod_{i\in \A(z)}
\theta_i(y')\equiv \iota_{y'}(\alpha)
\prod_{i\in \A(z)}\theta_i(y) \mod \theta_\ell(y)
\]
by Lemma~\ref{lemm:5}, and 
\[
\sgn(y)+\ly\sgn(y')=0
\]
by (\ref{eqn:sign of b}), so 
the numerator of the right hand side of the identity (\ref{eqn:y,y'}) 
is divisible by $\theta_\ell(y)=cf_1$. 
This means that we can arrange the left hand side of 
(\ref{eqn:y,y'}) with a common denominator in which $f_1$ does not appear 
as a factor. 
Since elements in $Q^c$ appear pairwise like this, one has 
\begin{equation*} \label{eqn:Qc}
\sum_{y\in Q^c}\frac{\sgn(y)\iota_y(\alpha)}{m(y)\prod_{i\in A(y)}
\theta_i(y)}=\frac{g_2}{\prod_{j=2}^N f_j}
\end{equation*}
with $g_2\in \Q[\Theta]$. 
This together with (\ref{eqn:Q}) implies that the numerator $g$ in 
(\ref{eqn:common}) is divisible by $f_1$. 
\end{proof}

Since $\iota_y$ is a $\Q[\Theta]$-module map, so is $\indT$.  Therefore 
$$\indT\colon A_P\to \Q[\Theta]$$ 
induces a homomorphism 
\begin{equation} \label{eqn:ind}
\ind \colon A_P/(\Theta)\to \Q.
\end{equation}
This map decreases degrees by $d$ because $\indT$ does.

\bigskip
\section{Gorenstein* simplicial posets} \label{sect:gore}

We shall prove Theorem~\ref{theo:main} in this section. 
Let $\k$ be an arbitrary field.  Suppose that a simplicial poset $P$ is 
Gorenstein* over $\k$, i.e., the order complex $\Delta(\oP)$ 
of $\oP=P-\{\hat 0\}$, which is a simplicial complex, 
is Gorenstein* over $\k$. 
According to Theorem II.5.1 in 
\cite{stan96}, a simplicial complex $\Delta$ of dimension $d-1$ 
is Gorenstein* over $\k$ 
if and only if for all $p\in |\Delta|$, 
\[
\widetilde H_q(|\Delta|,\k)\cong H_q(|\Delta|,|\Delta|-p;\k)
\cong \begin{cases} 
\k,\quad&\text{$q=d-1$,}\\
0,\quad&\text{$q<d-1$.}
\end{cases}
\]
Therefore, it follows from the universal coefficient theorem 
(\cite[Corollary 55.2]{munk84}) 
that if a simplicial poset $P$ is Gorenstein* over $\k$, then it is 
Gorenstein* over $\Q$.  In the sequel we may assume $\k=\Q$. 
According to Theorem II.5.1 in \cite{stan96} again, 
$\Delta(\oP)$ is an orientable pseudomanifold, 
so the assumption in Section~\ref{sect:index} 
is satisfied for the Gorenstein* simplicial poset $P$ 
because $\Delta(\oP)$ is the barycentric subdivision of the 
CW-complex $\Gamma(P)$. 

Since a Gorenstein* simplicial poset is Cohen-Macaulay, 
$h_i=h_i(P)$ agrees with the dimension of the homogeneous part of 
degree $i$ in $A_P/(\Theta)$, see 
the proof of Theorem 3.10 in \cite{stan91}.  
Therefore, if $h_i=0$ for some $i$ 
$(1\le i\le d-1)$, then a product of $d$ elements in $(A_P)_1$ vanishes 
in $A_P/(\Theta)$, 
in particular, the product is zero when evaluated by the index map 
in (\ref{eqn:ind}).  

We take a subset $I$ of $\{1,\dots,n\}$ with cardinality $d$ such that 
$I=\A(y)$ for some $y\in P_d$.  
If $\A(y)=\A(y')(=I)$, then $m(y)=m(y')$ by (\ref{eqn:m(y)=m(y')}). 
Therefore we may write $m(y)$ as $m_I$.  Since 
\[
\iota_y(\prod_{i\in I}x_i)=\begin{cases} \prod_{i\in \A(y)}\theta_i(y)
\qquad &\text{if $\A(y)=I$,}\\
0 &\text{otherwise}
\end{cases}
\]
by (\ref{eqn:iota}) and (\ref{eqn:theta}), we have 
\[
\indT(m_I\prod_{i\in I}x_i)=\sum_{\A(y)=I}\sgn(y)\in \Q 
\]
by (\ref{eqn:index}).  Hence, if we regard 
$m_I\prod_{i\in I}x_i$ as an element in $A_P/(\Theta)$, 
then we have 
\begin{equation} \label{eqn:final}
\ind(m_I\prod_{i\in I}x_i)=\sum_{\A(y)=I}\sgn(y).
\end{equation}

Now suppose that $h_i=0$ for some $i$ $(1\le i\le d-1)$. 
Then the left hand side 
of (\ref{eqn:final}) is zero as remarked above.  This means that 
(since $\sgn(y)=\pm 1$) there must be an 
even number of elements $y\in P_d$ with  
$\A(y)=I$ at the right hand side of (\ref{eqn:final}). 
Since $I$ is arbitrary, we conclude that 
$f_{d-1}$(the number of elements in $P_d$) is even. 
This together with (\ref{eqn:fn-1}) 
completes the proof of Theorem~\ref{theo:main}. 

\begin{rema} 
An element corresponding to $\w_T$ in Section~\ref{sect:topo} 
is $\sum_I m_I\prod_{i\in I}x_i$, where $I$ runs over all subsets of 
$\{1,\dots,n\}$ with cardinality $d$ and $m_I$ is understood to be zero 
if there is no $y\in P_d$ such that $I=\A(y)$. 
\end{rema}

\providecommand{\bysame}{\leavevmode\hbox to3em{\hrulefill}\thinspace}


\begin{thebibliography}{GZ}


\bibitem{duva94}
A.~M.~Duval, \emph{A combinatorial decomposition of simplicial complexes}, 
Israel J. Math. {\bf 87} (1994), 77--87.

\bibitem{fult93}
W.~Fulton, 
\emph{An Introduction to Toric Varieties}, 
Ann. of Math. Studies, {\bf 113}, Princeton Univ. Press, Princeton, NJ, 
1993. 

\bibitem{GZ}
V.~Guillemin and C.~Zara, \emph{Equivariant de Rham theory and graphs}, 
Asian J. Math. {\bf 3} (1999), 49--76. 

\bibitem{HM}
A.~Hattori and M.~Masuda, \emph{Theory of multi-fans}, Osaka J. Math. 
{\bf 40} (2003), 1--68.

\bibitem{ma-pa03}
M.~Masuda and T.~E.~Panov, \emph{On the cohomology of torus manifolds}, 
preprint. 

\bibitem{munk84}
J.~R.~Munkres, \emph{Elements of Algebraic Topology}, 
Addison-Wesley Publishing Company, 1984. 

\bibitem{stan80}
R.~P.~Stanley \emph{The number of faces of a simplicial convex polytope}, 
Adv. in Math. {\bf 35} (1980), 236--238.


\bibitem{stan91}
R.~P.~Stanley,
\emph{$f$-vectors and $h$-vectors of simplicial posets},
J. Pure Appl. Algebra {\bf 71} (1991), 319--331.

\bibitem{stan96}
R.~P.~Stanley,
\emph{Combinatorics and Commutative Algebra}, second edition,
Progress in Math. {\bf 41}, Birkh\"auser, Boston, 1996.


\end{thebibliography}
\end{document}